# *Palindromes In Sturmian Strings*

Ayşe Karaman
makalelr@yahoo.ca

**Abstract**

Let $p$ be a maximal palindrome in a Sturmian word $s=ul_1pl_2v$ so that $p$ is a palindrome and $l_1pl_2$ is not for letters $l_1$ and $l_2$. Let $\alpha_{(p,p')}$ be a morphism mapping letters $a$ and $b$ respectively to $a^p b$ and $a^{p'}b$, $|p-p'|=1$. In this paper, we characterize the palindromes in a Sturmian word and show that the number of maximal palindromes in a Sturmian word $X=\alpha_{(p,p')}(Y)$ for finite $Y$ and thus $X$ is $2|X|-2|Y|$. We show that the set of maximal palindromes in a finite Sturmian word $X$ has the cardinality $\sum_{i=1..n} \max(p_i, p'_i)$ where $X$ is characterized by subsequent mappings of $\alpha_{(p_i, p'_i)}$, $i=1..n$.

## 1. Introduction

Sturmian words are infinite words over a 2-letter alphabet which have the minimum possible factor complexity without falling into a repetition of a suffix infinitely many times. See for example [AS03] and [L02] for reference to the unique structural properties of Sturmian words. [B02]'s view describes Sturmian words as an arithmetic construction providing a bridge between combinatorics and number theory. A palindrome is a finite word that "reads" the same as its reverse. In this paper, we characterize the maximal palindromes in Sturmian words. We enumerate the distinct maximal palindromes in a given Sturmian word and their occurrences in the word. We also show an iteration which we call $\wp$ to generate the maximal palindromes, and in the limit the bistandard bi-infinite word having the same factors as a given Sturmian word. The next section includes the preliminaries and a thorough review of literature on palindromes in Sturmian words. Section 3 presents our results. In Section 4 we outline the morphic approach used throughout this study, and a linear algorithm to compute the maximal palindromes in a Sturmian word before we conclude.

## 2. Preliminaries

For a finite word $X$, $|X|$ and $|X|_l$ denote respectively the length of $X$ and the number of the occurrences of $l$ in $X$ where $l$ is a letter in the alphabet. $|X|_l$ is specifically named as *weight* of $X$ when $l=b$. We use the notation $X[i]$ to represent the letter at position $i$ of word $X$. Similarly, $X[i..j]$ denotes the *factor* $Y$ of $X$ where $|Y|=j-i+1$ and $Y[k] = X[i+k-1]$ for all $k=1..j-i+1$. A factor Y of X is named as a *proper factor* when $|Y|<|X|$. The *reverse* of a finite word $X$ is the word $X^R$ where $X[i] = X[|X|-i+1]$ for all $i=1..|X|$. A finite word $X$ is said to be a *palindrome* when $X=X^R$. Fact($s$) denotes the set of finite factors of $s$.

A Sturmian word is a right-infinite word which has the minimum factor complexity and which is aperiodic. Minimum factor complexity is attained on 2 letters among the aperiodic words and Sturmian words are defined on a binary alphabet. In this text, we use $A=\{a,b\}$ for the binary alphabet. Let $\alpha_{(p,p')}$ for a *valid parameter pair* $(p,p')$ satisfying $|p-p'|=1$, $p,p'\geq 1$ define the morphism on a 2-letter alphabet to map each letter $a$ to the sequence $a^p b$, and each letter $b$ to the sequence $a^{p'}b$. We equivalently use $\alpha$ to refer to $\alpha_{(p,p')}$ with any valid parameter pair $(p,p')$ whenever there is no need to be specific on the value of this pair. The morphism α is a Sturmian morphism mapping a Sturmian word to another [BS93]. Let $\pi=\{(p_1, p'_1), (p_2, p'_2), ... \}$ with each $(p_i, p'_i)$ a valid parameter pair denote a *defining sequence*. Let also $\alpha_\pi = ... \bullet \alpha_i ... \bullet \alpha_2 \bullet \alpha_1$ for some defining sequence $\pi$ where $\alpha_i = \alpha_{(p_i, p'_i)}$. In [K98] we present a combinatorial proof verifying that α is a Sturmian morphism, and show that every Sturmian word is a factor of $\alpha_\pi(a)$ for some defining sequence $\pi$. This characterization allows a vertical view of Sturmian words as a composition of levels thru which the structural properties are transferred. In [K10], we used this view to compute the runs in a Sturmian word. In this paper, we use the same approach for the computation of palindromes so that we trace the palindromes as they first appear in their original levels, then we compute how they reflect at further levels and finally the ultimate level of the word. On a given a finite defining sequence $\pi$ of size $|\pi|=n$, we say that a word $S_i$ is of *level i* if $S_i = \alpha_\pi(s_i) = \alpha_i(\alpha_{i-1}( ... (\alpha_1(s_i)) ...))$ for some $i \leq n$ and balanced word $s_i$. In this case, we call level $n$ the *ultimate level*. Observe that a Sturmian word $S_i$ of level $i$ is a concatenation of blocks $a^{p_i} b$ and $a^{p'_i} b$, which we call *full blocks*. A factor in $S_i$ which is a proper factor of a full block is called a *partial block*. In $S_i$, let $p_{min} = \min\{ p_i, p'_i \}$ and correspondingly $p_{max} = p_{min}+1$. We name the full blocks $a^{p_{min}} b$ and $a^{p_{max}} b$ respectively as *short* and *long* blocks. Observe in this case that $a^{p_{min}} b$ has both full and partial block occurrences in $S_i$. In $S_i$, $\alpha_i(a)$ and $\alpha_i(b)$ are respectively the *repeating* and *non-repeating* blocks. We name as *a-sequence* a factor $X$ of a Sturmian word $S=UXV$ where $|X|_b=0$, $U[|U|]=b$ and $V[1] = b$. According to this, the *a*-sequences of a Sturmian word $S=\alpha_{(p,p')}(S')$ are either of the form $a^{p_{max}}$ or $a^{p_{min}}$. We define $\alpha^{-1}$ as the inverse of α on *block-complete* words, i.e., words that are concatenation of full blocks so that $\alpha^{-1}(\alpha(X))=\alpha(\alpha^{-1}(X))$. According to this, $\alpha^{-1}(X)=Y$ iff $\alpha(Y)=X$.

The following property is directly implied by α and useful for the "arithmetic" through the levels imposed by this morphism:

**Property 1:** Let $X= \alpha_{(p,p')}(Y)$ where $Y$ and thus $X$ are finite Sturmian words. Then,
i.) $|X|=(p+1)|Y|_a + (p'+1)|Y|_b$,
ii.) $|X|_a=p|Y|_a + p'|Y|_b$,
iii.) $|X|_b=|Y|$. □

Given a palindrome $X$, the *center* of palindrome is the letter at position $(|X|+1)/2$ if $|X|$ is odd, and the 2-letter sequence at consecutive positions $|X|/2$ and $(|X|+2)/2$ otherwise. Observe that a 2-letter-center is always *aa* while a single letter in the center is either *a* or *b*. We commonly call the centers with the letters involved as *a-center*, *b-center* or *aa-center*. By definition, every center is also a palindrome. Let $P_1$ and $P_2$ be two palindromes so that $P_2$ is a subtring of $P_1$. Furhermore, $P_1$ and $P_2$ share a center in their occurrences. We call such palindromes as *same-centric*. Given a palindrome $X$, each letter pair occurring at positions $i$ and $j$ of $X$ where $i+j=|X|+1$ is called *twin letters*. Similarly, we call two factors $U$ and $V$ of a palindrome $X=W_1UW_2VW_3$ *twin sequences* whenever $|U|=|V|$ and $|W_1|=|W_3|$, in which case $U$ and $V$ are reverse words in same "distance" from the palindrome center.

A palindrome $P$ is *maximal in a word S* if $P$ has an occurrence in $S=UPV$ so that $U[|U|]PV[1]$ is not a palindrome. We refer by *non-maximal occurrence* of a maximal word $P$ in $S$ to an occurrence of $P$ in $S$ as a factor of $lPl$ for some $l \in \{a,b\}$.

The set PER, equivalently the set of *central words* in the Sturmian literature refers to those set of palindromes which have maximal occurrences in some Sturmian word. We denote by $M_S$ the set of palindromes that are maximal in $s$. $M_{(S,c)}$ denotes the set of maximal palindromes with center $c$ in $S$. $K_{(S,c)}$ represents the set of maximal palindrome occurrences with center $c$ in $S$, i.e., $K_{(S,c)}=\{(i,c) \mid c$ is a palindrome center in $S$ so that $S[i] = c[1]\}$. $K_S$ is the set of all maximal occurrences in $S$. Specifically, $K_S = K_{(S,a)} \cup K_{(S,aa)} \cup K_{(S,b)}$ and $M_S = M_{(S,a)} \cup M_{(S,aa)} \cup M_{(S,b)}$. Note the distinction between $M_S$ and $K_S$ that $M_S$ is the set of *distinct* maximal palindromes in $S$ while $K_S$ shows their occurrences in the word. Note the distinction between the elements of PER and maximal words in that the set PER characterizes a set of palindromes across the class of Sturmian words while the set of maximal palindromes is attributed to some specific word in the Sturmian class. Also note that not every factor of a Sturmian word that is also a member of PER is maximal in that word, i.e., $M_s$ is a proper subset of PER∩Fact($s$), as exemplified by $a^{P_{max}}$ ∈PER∩Fact($s$) \ $M_s$ whenever $a^{P_{max}}b$ is the long block in $s$.

## 2.1. Literature Survey

This section reviews the literature on the palindromic characteristics of Sturmian words and shows the place of our results on Sturmian palindromes with respect to those in literature.

*Structural properties – set PER, maximal palindromes and singular words*: The set PER directly connects to the set of finite standard words and well-characterizes the Sturmian words to dominate the view on the palindromes and even on the structure of the class of Sturmian words. Set PER is initially hinted by Coven and Hedlund [CH73] with some of its properties noted below. de Luca and Mignosi [DM94] fully defined this set and characterized many of its properties as a special set of balanced sequences. The following statements are equivalent and review the literature on this perspective. By definition, $a^k$∈PER for all $k≥1$. In all of the below cases, the domain of PER is restricted without loosing generality to the words with at least one letter of each of the binary alphabet.
- $p$∈PER.
- $pl_1l_2=qr \Leftrightarrow |q|+2$ and $|r|$-2 are coprimes where $p$, $q$ and $r$ form a unique set of palindromes on $\{a,b\}^*$ with $|q|>|r|$, $l_1,l_2\in\{a,b\}$, $l_1\neq l_2$ [P88] (cf. [BD97]). In this case, $s=pl_1l_2$ for some standard word $s$ and all finite standard words are of this form [DM94]. Any circular shift of $s$ is either a palindrome or a concatenation of two palindromes [CW03].
- $p=u^m$ and $p=v^n$ for exactly two factors $u\neq v$ where $|u|$ and $|v|$ are coprimes, $|p|=|u|+|v|$-2 and such decomposition is unique by [DM94] and direct implication of [CH73]. In this case, $u= q_1l_1l_2$ and $v= q_2l_2l_1$ for $q_1$, $q_2$∈PER and unequal letters $l_1$ and $l_2$ in the alphabet. Let $w^{(-)}$ denote the smallest palindrome having $w$ as the suffix. Then, $p=(lw)^{(-)}$ for $w$∈PER, $l$∈$\{a,b\}$. Furthermore, $(ap)^{(-)},(bp)^{(-)}$∈PER to define a constructive description of this set [D97].

- *p* is strictly bispecial, *i.e.*, *apa, bpb, apb* and *bpa* are Sturmian factors, equivalently finite balanced words [DM94]. *apb*, *bpa* and *lpl* for *l=a* or *l=b* are three factors of some Sturmian word *s* by direct implication of [CH73].

Singular words constitute a specific class of palindromes associated with the finite standard words and have been used to analyze the structural properties of Sturmian words. Let *s* be the standard word of some slope $\theta$. Let $s_i$ denote the finite standard words of *s*. As stated above, every $s_i = pl_1l_2$ is of the form where $p \in$ PER and by implication $p \in M_s$, $l_1$ and $l_2$ are unequal letters in the alphabet. The set of singular words associated with *s* or any of the Sturmian words of slope $\theta$ is exactly the set of palindromes of the form $l_1pl_1$ [CW03]. Singular words have been studied in [M99, CW03] to analyze the structural properties of Sturmian words. Intuitively, singular words $l_1pl_1$ as above mark the non-maximal occurrences of maximal word *p* in *s* and reverse the properties of these words. [M99] showed that a singular word is never a concatenation of two palindromes while this is not a general property of the words of PER. [CW03] showed that no circular shift of a singular word is balanced—note to the contrary that all circular shifts of a word of PER are balanced. [G06] extends the comprehensive results in [CW03, M99] and gives explicit formula to return the positions of a given palindrome in a standard word. The following decompositions are shown to characterize the standard Sturmian words by the use of the class of singular words.

- A standard word *s* is of the form $s = q_1q_2q_3...$ where $q_i = lu$, with *u* the maximum proper prefix of $s_{i-1}^{a_i-1} s_{i-2}$ and *l* the letter unequal to the last letter of the standard word $s_i = s_{i-1}^{a_{i+1}} s_{i-2}$ [M99, verified by CW03]. Observe that $q_i$ are exactly the singular words of *s* itself only if *s* is the Fibonacci word with partial coefficient values of 1 [WW94].
- A standard word is a composition of three palindromes for infinitely many such palindrome triples. Specifically, the consecutive occurrences of a singular word $w_n$ in a Sturmian word *s* are separated by $w_{n+1}$ and $w'_n$ where $w_{n+1}$ *and* $w'_n$ are singular words with $w_{n+1}$ singular in *s* and $w'_n$ singular in some $s' \neq s$ as precisely described in [CW03].

*Number of palindromes – palindromic complexity and palindromic richness*: Before going into this topic, note the significant property known by [CH73] and verified by [M89] that all Sturmian words are closed under reversal, i.e., a finite word *w* is a factor of a Sturmian word *s* iff the reverse of *w* is a factor of *s*. *Palindromic complexity*, analogous to the factor complexity, measures the number of distinct palindromes of a given length in a word whereas *palindromic richness* refers to the property that a word has maximum number of palindrome factors possible. Specifically, a word is called palindromically rich iff every finite factor *u* has |*u*| distinct number of non-empty palindromes which is the maximum possible. Equivalently [GJWZ09], a word is called palindromically rich iff every factor having two consecutive occurrences of a given palindrome as a border is itself a palindrome. Palindromic richness of a word is implied by its being closed under reversal [BDGZ09] and all Sturmian words are rich [DJP01] (also note the direct implication by the closure-under-reversal property). The cardinality of palindromes in a word is also approached from the view of palindrome complexity. [ABCD03] shows an upper bound for palindrome complexity in linear terms of the factor complexity of a word in the general class of aperiodic infinite words. [BMP07] studies the problem on uniformly recurrent words and shows an upper bound for palindrome complexity of the words closed under reversal. [BDGZ09] uses the bound in [BMP07] to imply palindromic richness when attained in any class of words closed under reversal which includes the class of Sturmian words. It is known by [DP99] that the class of Sturmian words are exactly the words that have exactly one (two) palindromic factor of each even (odd) length. In [DP99], the authors also question and prove the existence of a morphism to generate the palindromes of a given Sturmian word. [DM94] shows for a given length $n$ that the number of palindromes of length $n$ in PER is $\phi(n+2)$ where $n$ is Euler totient function. [DD06] extends this count to the general domain of balanced palindromes of a given length.

*Palindrome prefixes and palindrome density*: The *palindromic prefixes* of a standard word $s$ are exactly the set of palindromes which have a maximal occurrence in $s$, and each palindromic factor of $s$ is a same-centric factor of such a prefix of $s$ [D97]. It is known by [BR05] that palindromic factors of a Sturmian word and thus such prefixes of the standard word are arbitrarily long. [F06] uses palindrome prefixes to define a measure of the palindrome density of a word. This definition puts the infinite words into a scale of which the words with finitely many palindrome prefixes and the purely periodic words are at respectively the low and high extremes. According to this measure, the Fibonacci word has the highest palindrome density among the general class of aperiodic words. Furthermore, for each irrational $\theta$, the standard word has the highest (defined) palindrome density among the words of the same slope.

In this paper, a formula is given for the number of *distinct* maximal palindromes across their varying palindrome centers in a given Sturmian word. Furthermore, it is shown that the number of *occurrences* of maximal palindromes in a Sturmian word $s = \alpha_{(p,p')}(s')$ is $2(|s|-|s'|)$ and thus is linear in the string length. The results of this paper relate to those in literature as follows:
- [DM94] and this work study the number of elements of PER which respectively are any balanced sequences of a given length, and which have maximal occurrences in a given Sturmian word. The two papers study the Sturmian factors respectively in the domain of entire balanced sequences and in the factor set of a given Sturmian word. [DM94] and this paper provide solutions to two different problems and, to the best of our knowledge, the implications of either of them on the other is not studied in literature. [DD06] further extends the results in [DM94] to count the number of any balanced palindromic sequences of a given length without restricting to distinct Sturmian sequences as this paper does and to the members of set PER as in [DM94] and in this paper.
- [G06] enumerates the positions of a given palindrome in a standard Sturmian word given by its slope. This paper describes a linear time algorithm to compute the sequences and positions of the maximal palindromes in a Sturmian word given by its defining sequence. The results in [G06] apply with direct implication to the set of maximal palindromes of the word. The algorithm described in this work additionally enumerates the maximal palindromes in the Sturmian word it processes on.

- Our results verify in process the following:
    - Every maximal palindrome also occurs as a non-maximal palindrome in a Sturmian word [DM94],
    - Palindromic factors of a Sturmian word are arbitrarily long [BR05],
    - Sturmian words are closed under reversal [CH73, M89].

## 3. Findings

As we stated Section 2, the formulation of a Sturmian word as an image of parametric morphisms $\alpha_{(p,p')}(Y)$ allows the view of such a word as a composition of levels each identified by the associated valid parameter pair $(p,p')$. This approach facilitates the characterization of structural patterns in Sturmian words so that an occurrence of a pattern is first identified at the level it originally appears, and the reflection of this original form is computed in the further levels and eventually the ultimate level.

We first establish the link between the levels to see the how a palindrome $X$ reflects in $\alpha(X)$ – the next level.

**Lemma 1:** Let $X$ be a palindrome in a Sturmian word. Then $b\alpha(X)$ is a palindrome. Conversely, if $Xb$ is block-complete with $X$ a palindrome, then $\alpha^{-1}(Xb)$ is a palindrome.

**Proof:** Let $X$ be a factor of a Sturmian word $S$. Observe that $\alpha(X)$ is always preceeded by a letter $b$ in $\alpha(S)$. (Note the case that $X$ is a prefix of $S$ and $\alpha(X)$ is a prefix of $\alpha(S)$— thus this condition does not hold without contradicting the generality of this proof.) Consider every $a$-sequence twin pairs $a^{m1}$ and $a^{m2}$ preceeding and succeeding respectively the center of $b\alpha(X)$. For every such pair, $m1 = m2$ since $a^{m1}b$ and $a^{m2}b$ are the mappings of the twin-letters in $X$.

Conversely, every twin $a$-sequence $a^{m1}$ and $a^{m2}$ in a palindrome $bX$ are necessarily of the same length. Considering $Y=\alpha^{-1}(X)$, $a^{m1}b$ and $a^{m2}b$ correspond to twin letters in $Y$ for every such $a^{m1}b$ and $a^{m2}b$ sequence pair in $bX$. □

**Corollary 1.a:** Let $c$ be the center of a palindrome $X$ in a Sturmian word. Then the center of $b\alpha_{(p,p')}(X)$ is the center of $b\alpha_{(p,p')}(c)$.

**Proof:** The center of a palindrome is a palindrome by definition. Follows from Lemma 1. □

We say that a palindrome with center $c \in \{a,b,aa\}$ is *original* in a Sturmian word $X=\alpha_{(p,p')}(t_1 l t_2)=x_1 c x_2$ unless $c$ is the center of a palindrome $b\alpha_{(p,p')}(l)$ while $\alpha_{(p,p')}(l) = ucv$, $x_1 = \alpha_{(p,p')}(t_1)u$, $x_2 = v\alpha_{(p,p')}(t_2)$ and $l \in \{a, b, aa\}$. A palindrome is called a *reflection* if it is not original at the current level of the word. We use the notation $O_X$ to denote the set of original palindromes in $X$. We further establish in this text that a maximal palindrome $X$ in a Sturmian word is original iff $|X|_b \leq 1$.

By definition, each letter and each $aa$ sequence is a palindrome-center defining exactly one instance of a maximal palindrome. Note here this directly justifies the fact that *the number of*

*maximal palindromes in a Sturmian word is linear in the word length.* The concern is whether a center is original or the center of a reflection of a palindrome from an earlier level. The following two proofs are trivial and tell this distinction.

**Corollary 1.b:** Let $c \neq b$ be the center of a palindrome $x$ in a Sturmian word. Then $x$ is original in its level iff $c$ is not the center of its $a$-sequence. □

**Corollary 1.c:** A center $b$ is the first letter $b$ in a given $\alpha(l_1 l_2)$ where $l_1, l_2 \in \{a,b\}$ and $l_1 \neq l_2$ iff $b$ is an original center. □

Lemma 2 characterizes the reflection of a maximal palindrome in the next level. Before this proof, note the following property that a maximal palindrome $X$ with $|X|_b > 0$ in a Sturmian word $S$ has $a^{p_{min}}b$ and $ba^{p_{min}}$ respectively as a prefix and suffix for $aXb$ (equivalently $bXa$) is a factor of $S$. However, the opposite of this is not true—not every factor of the form $Uba^{p_{min}} = a^{p_{min}}bV$ of $S$ is maximal in $S$.

**Property 2:** $X \in M_S$ and $|X|_b > 0 \Rightarrow X = Uba^{p_{min}} = a^{p_{min}}bV$ where $S$ is a Sturmian word and is the short block in $S$.

**Proof:** $X \in M_S$ implies that $aXb, bXa \in Fact(S)$. $aX, bX \in Fact(S)$ and $|X|_b > 0$ imply that $a^{p_{min}}b$ is a prefix of $X$ since any other prefix $a^k b$ for $k \neq p_{min}$ resulted $ba^k b$ or $ba^k b$ be unbalanced with some factor of $S$. With a similar argument, $X$ has the suffix $ba^{p_{min}}$. □

**Lemma 2:** A palindrome $X$ is maximal in a Sturmian word $S$ iff the palindrome $a^{p_{min}}b\alpha_{(p,p')}(X)a^{p_{min}}$ for $|p-p'|=1$ and $p_{min} = \min\{p,p'\}$ is maximal in $\alpha_{(p,p')}(S)$.

**Proof:** Let $X$ be the factor of $l_3 l_1 X l_2$ where $l_1, l_2, l_3 \in \{a,b\}$. By condition, $l_1 \neq l_2$. By Lemma 1, $b\alpha_{(p,p')}(X)$ is a palindrome and so is $Wb\alpha_{(p,p')}(X)W^R$ for any finite word $W$.

Consider $\alpha_{(p,p')}(l_3 l_1 X l_2) = A_3 b A_1 b \alpha_{(p,p')}(X) A_2 b$, where $A_1, A_2$ and $A_3$ are $a$-sequences. Since $l_1 \neq l_2$, $A_1 \neq A_2$ where $A_1, A_2 \in \{a^{p_{min}}, a^{p_{max}}\}$. Thus, $a^{p_{min}}$ in the suffix of $A_1$ and prefix of $A_2$ is a twin sequence pair extending the palindrome $b\alpha_{(p,p')}(X)$. And, the letters $A_1[1]$ and $A_2[|A_2|]$ in twin positions are different letters since $A_1 \neq A_2$.

Conversely, let $X'$, where $|X'|_b \geq 1$ be a maximal palindrome that is a reflection in $\alpha_{(p,p')}(S)$. Then, by Property 2, $X'$ has the prefix $a^{p_{min}}b$ and suffix $ba^{p_{min}}$. Let $X' = a^{p_{min}}bYa^{p_{min}}$. $Y$ is either empty or block-complete and thus $\alpha^{-1}(Y)$, is defined. $a^{p_{min}}ba^{p_{min}}$ is by definition original in the current level of $S$. When $Y$ is not an empty word, the twin $a$-sequences in $Y$ correspond to twin letters in $\alpha^{-1}(Y)$, which makes $\alpha^{-1}(Y)$ a palindrome. To show $\alpha^{-1}(Y)$ is a maximal palindrome in $S$, examine the full-blocks $f_1$ and $f_2$ respectively preceeding and succeeding $Y$ in $\alpha(S)$. $f_1 \neq f_2$, otherwise $X' = a^{p_{min}}bYa^{p_{min}}$ is a factor of $la^{p_{min}}bYa^{p_{min}}l$ for $l \in A$ and $X'$ is not maximal. The unequal full blocks surrounding $Y$ are images of unequal letters surrounding $\alpha^{-1}(Y)$, concluding the proof. □

By Lemma 2, the parametric mapping on words $\wp_{(p,p')}(X) = a^{p_{min}} b \alpha_{(p,p')}(X) a^{p_{min}}$ where $p_{min} = \min\{p,p'\}$ is a palindrome-function mapping a maximal palindrome in a Sturmian word $S$ to another in $\alpha_{(p,p')}(S)$ through the level defined by $(p,p')$. Recall the definition of a maximal palindrome $X$ that $aXb$ (and thus $bXa$ from reverse-closure property) are factors of a given word. Note also the definition of bispecial words of a given word $S$ that $w$ is called *bispecial* in $S$ if $aw$, $bw$, $wa$ and $wb$ are factors of $S$, and thus the maximal words of $S$ are bispecial in $S$ [DM94]. In our context, we equivalently name the bispecial words of a Sturmian word $S$ as *bistandard* to refer to the analogy between the finite and infinite standard words of the right-infinite class of Sturmian words, and the such words of the bi-infinite case. Lemma 2 above also results in that $\wp$ is a bistandard function to map bistandard words into others. The limit of the word $\wp^i$ as $i$ goes to infinity then produces a bistandard infinite word.

**Theorem 1:** Let $X = \alpha_{(p,p')}(Y)$ be a finite Sturmian word. Then, the following are correct:

i.) $|K_{(X,a)}| = p|Y|_a + p'|Y|_b$.
ii.) $|K_{(X,b)}| = |Y|$,
iii.) $|K_{(X,aa)}| = (p-1)|Y|_a + (p'-1)|Y|_b$,

iv.) $|K_X| = 2|X| - 2|Y| = 2|Y|_a + 2p'|Y|$ when $p > p'$, and $K = 2|Y|_b + 2p|Y|$ when $p < p'$.

**Proof:**
Property 1 directly justifies cases (i) and (ii).

Clearly, the number of distinct $aa$ sequences in an $a$-sequence of length $|A|$ is $|A|-1$. In $x$, the $a$-sequences are either of length $p$ or $p'$, respectively being the reflections of letters $a$ and $b$ in $y$. So, $|Y|_a$ and $|Y|_b$ are the number of $a$-sequences respectively of length $p$ and $p'$, and the number of $aa$-centers in $x$ is $(p-1)|Y|_a + (p'-1)|Y|_b$, justifying case (iii).

The number of maximal palindrome occurrences is the sum of the number of palindromes of all possible centers, which from above is:

$|K_X| = (p|Y|_a + p'|Y|_b) + (|Y|) + [(p-1)|Y|_a + (p'-1)|Y|_b] = 2(p|Y|_a + p'|Y|_b)$.

When $p > p'$, $|K_X| = 2|Y|_a + 2p'|Y|$.
When $p < p'$, $|K_X| = 2|Y|_b + 2p|Y|$.

$2|X| = 2(p|Y|_a + p'|Y|_b) + 2(|Y|_a + |Y|_b)$, and $|K_X| = 2|X| - 2|Y|$. □

**Lemma 3:** Let $X$ be a finite Sturmian word so that $X = \alpha_{(p_x,p'_x)}(Y) = \alpha_{(p_x,p'_x)}(\alpha_{(p_y,p'_y)}(Z))$. The total number of maximal original palindrome occurrences in $X$ is

$2(p_x-1)(|Z|_a p_y + |Z|_b p'_y) + 2|Z|p'_x$

**Proof:**
The number of original palindromes in $X$ is:

$$\begin{aligned}
C_X &= |K_X| - |K_Y| \\
&= 2|X|-4|Y|+2|Z| \quad \text{// by Theorem 1} \\
&= 2(p_x+1)|Y|_a + 2(p'_x+1)|Y|_b - 4(p_y+1)|Z|_a - 4(p'_y+1)|Z|_b + 2|Z|_a + 2|Z|_b \\
&= 2(p_x-1)(|Z|_a p_y + |Z|_b p'_y) + 2|Z|p'_x \quad \text{// by Property 1 and simplification} \quad \square
\end{aligned}$$

The above two proofs enumerate the occurrences of maximal palindromes in a Sturmian word. The following proofs characterize the set of maximal palindromes in a given word.

**Lemma 4:** Let $P$ be a maximal palindrome in a Sturmian word $S$. $P$ is original in its level iff $|P|_b \leq 1$.

**Proof:**
Let $P$ be a maximal original palindrome where $|P|_b > 1$. If $|P|_b$ is even, the center $a$ or $aa$ of the a-sequence is not original by Corollary 1.b. If $|P|_b$ is odd and thus $P$ is a b-centered palindrome, then center $b$ is the first $b$ in $\alpha(l_3 l_4)$ where $l_3, l_4 \in \{a,b\}$ and $l_3 \neq l_4$ (Corollary 1.c). But when $l_3 \neq l_4$, $l_3$ and $l_4$ are mapped to 2 different a-sequences around center $b$—implying that $b$ is a center of a palindrome $P$ where $|P|_b = 1$, a contradiction. Palindrome $P$ satisfying $|P|_b \leq 1$ is of the form $a^n$ or $a^n b a^n$ in both of which cases it is original. $\square$

The findings ending by the results of Lemma 3 give us the number of occurrences of maximal palindromes. We use Lemma 4 in the following proof to provide the exact number of the distinct palindromic sequences regardless of their positions in a Sturmian word.

**Example:** Consider the standard Fibonacci word $f=abaababaabaababaababa...$ . The following are true for $f$:

  i.)   $\{(1,b), (5,b), (3,aa), (6,a)\} \subset K_f$.
  ii.)  $\{aba, a, abaaba\} \subset M_f$.
  iii.) $aabaa \notin M_f$ since $aaa$ is not a factor of $f$ and $aabaa$ appears only as a factor of $baabaab$ which always is the factor of a maximal palindrome and itself is not maximal. Similarly, $aa \notin M_f$ and thus only has non-maximal occurrences in $f$.
  iv.)  In $f$, letter $a$ is a maximal palindrome and has maximal occurrences whenever it is occurring as a factor of an $aa$ sequence. The rest of its occurrences in $bab$ sequences are as centers of reflections of palindromes from earlier levels. $\square$

**Lemma 5:** Let $X=\alpha_{(p,p')}(Y)$ be a Sturmian word. Then, the following are correct:

i.) $|M_{(X,b)} \cap O_X| = 1$
ii.) $|M_{(X,a)} \cap O_X| = (k+1)/2$ where $k \in \{p-1, p'-1\}$ and $k$ is odd
iii.) $|M_{(X,aa)} \cap O_X| = k/2$ where $k \in \{p-1, p'-1\}$ and $k$ is even

iv.) $|M_X \cap O_X| = \max(p,p')$.

**Proof:**
By Lemma 4, a maximal original palindrome in a Sturmian word is of weight 1 if it is $b$-centered, 0 otherwise.

Case (i) is trivial since a $b$-centered maximal palindrome is necessarily of the form $a^{p_{min}}ba^{p_{min}}$ and is a factor of $ba^{p_{max}}ba^{p_{min}}b$ or $ba^{p_{min}}ba^{p_{max}}b$. Note here that, $a^{p_{max}}ba^{p_{max}}$, although a palindrome, is never maximal since it only is the factor of $ba^{p_{max}}ba^{p_{max}}b$ in all its occurrences.

In an $a$-sequence, any proper prefix (or, equivalently, a proper suffix) is a maximal palindrome since preceeded and succeeded by different letters, and an original palindrome by Corollary 1.b since its center is not the center of the $a$-sequence itself. Since such proper prefix and suffix palindrome pairs in an $a$-sequence are equal, the elements of $M_{(X,a)} \cap O_X$ are every odd-length proper factor of any $a$-sequence of $X$. When $\max\{p,p'\}$ is even, those factors of length $\min\{p,p'\}$ of $a^{P_{max}}$ are original palindromes in their occurrences as a proper factor of an $a$-sequence of length $\max\{p,p'\}$, justifying case (ii) of the Lemma. Similar arithmetic leads to case (iii). Case (iv) is the sum of the 3 cases. □

The above Lemma states that the maximal original palindromes in a particular level of a Sturmian word are $a^{p_{min}}ba^{p_{min}}$ and $a^i$ for every $i=1..p_{min}$. The cardinality of $M_X$ is then trivial:

**Theorem 2:** Let $X = \alpha_{(p_1,p'_1)} \cdot \alpha_{(p_2,p'_2)} \cdot ... \cdot \alpha_{(p_n,p'_n)}(Y)$ be a Sturmian word. Then, $|M_X| = |M_Y| + \sum_{i=1..n} \max(p_i,p'_i)$. □

The following proof verifies a result by [DM94].

**Lemma 6:** Let $P$ be a maximal palindrome in a Sturmian word $X$. Then, the palindrome $P$ occurs in $X$ also as an non-maximal palindrome which is the proper factor of a same-centric maximal palindrome.

**Proof:**
A maximal palindrome of the form $a^{p_{min}}ba^{p_{min}}$ in Sturmian word $X$ has an occurrence in $X$ as a proper factor of $a^{p_{max}}ba^{p_{max}}$ when $(p,p')=(p_{max}, p_{min})$, and as a factor of $ba^{p_{min}}ba^{p_{min}}b$ when $(p,p')=(p_{min}, p_{max})$. $a^{p_{max}}ba^{p_{max}}$ and $ba^{p_{min}}ba^{p_{min}}b$ themselves are not maximal palindromes, being factors of $ba^{p_{max}}ba^{p_{max}}b$ and $aba^{p_{min}}ba^{p_{min}}ba$ respectively.

Likewise, a palindrome of the form $a^i$ has an occurrence as a factor of $a^{i+2}$ if $i+2 \leq p_{max}$, and of $ba^ib$ if $i=p_{max}-1$ – neither of which are maximal palindromes. Also, the morphism α is non-erasing and expands a word at every level, the reflection of a maximal palindrome further gaining a prefix extension and a suffix entension of respectively $a^{p_{min}}b$ and $a^{p_{min}}$, concluding for the proof. □

Corollary 6.a directly follows from Lemma 6 and justifies [BR05]'s result.

**Corollary 6.a.** A palindrome is arbitrarily long, i.e., for every positive integer $k$ and Sturmian word $S$, there is a palindrome $P$ in $S$ where $|P|>k$. □

**Corollary 6.b.** Any finite factor $X$ of a Sturmian word $S$ is the factor of a maximal palindrome in $S$. Furthermore, $X^R$ is a factor of $S$—justifying [CH73, M89].

**Proof:** Since every letter is a palindrome and extends in length (Lemma 2) and overlaps with the two palindromes on both sides in the next level (Lemma 6), every factor of a Sturmian word is a factor of a palindrome or a concatenation of palindromes in a given level—which in turn is a factor of one palindrome originated at some earlier level. □

## 4. Conclusion

The findings in this work used the α-morphic analysis [K98] of a Sturmian word given by its defining sequence π. This view takes a Sturmian word as a composition of levels each defined by a valid parameter pair $(p,p')$. In each such level, blocks act like letters so that $a^pb$ and $a^{p'}b$ are letters of one level mapped into blocks of the next. According to this view, the levels of a Sturmian word have identical properties that can be processed on so that once a pattern is identified at a level, its reflection can be computed iteratively at the further levels and finally in the ultimate level of the word. Thus, an algorithm to compute a pattern, which is a palindrome in the subject matter of this text (refer to [K10] as an application of this view on repetitions in Sturmian words) , is of two main steps:

1.) detect the pattern at the level it originates,
2.) compute how the pattern transmits through levels of α-morphic analysis on the defining sequence of the Sturmian word.

We generically name a pattern as "original" at the level it first exist as a result of step (1) above, and "reflection" at the further levels as results of step (2).

The cardinality of set of occurrences of maximal palindromes in a Sturmian word is linear in the length of the word. This can easily be checked by taking each letter and each *aa* sequence as the center of a maximal palindome in a Sturmian word. Using the results of Corollary 1.b and 1.c for step (1), a maximal palindrome occurrence can be detected by its position at the level it originates. The use of these results take constant time for each maximal palindrome. By the use of the results of Lemma 2 and Property 1 on the given defining sequence of a Sturmian word, the reflection of a specific palindrome at the ultimate level can directly be computed within a linear iteration on the Sturmian word. Note here that, locating the position of a specific palindrome at the ultimate level takes time linear to the in the size of the defining sequence of the word. However, each iteration through levels compute the reflection of all palindromes through the transmission of the level on the known composition of the Sturmian word-- maintaining the linearity of the overall computation time. And since the number of maximal palindrome occurrences are linear in a Sturmian word, an algorithm to compute each maximal palindrome occurrence by their poisitions in a Sturmian word has lienar complexity in the length of the word.

**Acknowledgements**

I thank my Dad for his pleasant company during the intellectual work into this paper.